\documentclass[10pt]{article}

\usepackage[latin1]{inputenc}
\usepackage{amsmath,amssymb,amsthm}
\usepackage{graphicx}
\usepackage{epsfig}

\newtheorem{theorem}{Theorem}[section]

\newtheorem{definition}[theorem]{Definition}

\newtheorem{remark}{Remark}

\newcounter{listagem}
\newcommand{\blista}{\begin{list}{\roman{listagem})}{\usecounter{listagem}}}
\newcommand{\elista}{\end{list}}

\newcommand{\beq}{\begin{equation}}
\newcommand{\eeq}{\end{equation}}
\newcommand{\beqn}{\begin{eqnarray}}
\newcommand{\eeqn}{\end{eqnarray}}

\def\le{<}
\def\ge{>}
\newcommand{\Cl}{{C \kern -0.1em \ell}}

\newcommand{\R}{\mathbb{R}}

\begin{document}

\title{Dirac and Laplace operators on some non-orientable conformally flat manifolds}
\author{Rolf S\"oren Krau{\ss}har \\ \\
{\small Arbeitsgruppe Algebra, Geometrie und Funktionalanalysis} \\
{\small Fachbereich Mathematik} \\ {\small Technische Universit\"at Darmstadt }\\{\small Schlo{\ss}gartenstr. 7}\\{\small 64289 Darmstadt, Germany.} \\ {\small E-mail: krausshar@mathematik.tu-darmstadt.de}} \maketitle


\begin{abstract}
In this paper we present an explicit construction for the fundamental solution to the Dirac and Laplace operator on some non-orientable conformally flat manifolds. We first treat a class of projective cylinders and tori where we can study monogenic sections with values in different pin bundles. Then we discuss the M\"obius strip, the Klein bottle and higher dimensional generalizations of them. We present integral representation formulas together with some elementary tools of harmonic analysis on these manifolds.
\end{abstract}


\section{Introduction}
The study of Dirac and Laplace operators on Riemannian manifolds has lead to a profound understanding of many geometric aspects related to these manifolds. In turn, Riemannian manifolds play a central role in several branches of modern physics. They appear in important  cosmological models, in general relativity theory, in the standard model of particle physics, in string theory and in general quantum field theory. The Dirac operator $D$ is a first order differential operator that acts on a vector bundle over such a Riemannian manifold. It linearizes the second order Laplacian, viz $D^2=-\Delta$. 
\par\medskip\par
Associated to the Dirac operator there is a rich function theory which allows us to tackle many boundary value problems of harmonic analysis on many subclasses of Riemannian manifolds. Consequently, it permits us to study the action of the Laplacian on these manifolds, too. One important tool in the analysis of the Dirac operator is the Cauchy integral formula. This formula allows us to calculate the values of functions in terms of given boundary data arising from practical measurements. Based on these representations important existence and uniqueness theorems for the solutions of boundary value problems on manifolds could be established, see for instance \cite{ca,cn}. 
\par\medskip\par
For a quantitative analysis one is particularly interested in explicit representation formulas for the appearing Cauchy kernels or harmonic Green's kernels.  However, this is a rather difficult task in general. In the recent years much success has been provided in the particular context of conformally flat spin manifolds. Conformally flat manifolds are Riemannian manifolds with a vanishing Weyl tensor. In turn, these are characterized by those manifolds which possess an atlas whose transition functions are conformal maps. Following the classical paper \cite{Kuiper}, some concrete examples can be constructed by factoring out a simply connected domain $G$ from $\mathbb{R}^n$ by a torsion free Kleinian group $\Gamma$ that acts totally discontinuously on $G$. The kernel then can be expressed in terms of automorphic forms on  $\Gamma$ that are in the kernel of the Dirac or Laplace operator, cf. \cite{KraHabil,KraRyan1,KraRyan2,BCKR}. In the case of an oriented conformally flat manifold one can consider one or several spinor bundles; we are dealing with a conformally flat spin manifold. Typical examples with these features are cylinders and tori, the Hopf manifold $S^1 \times S^{n-1}$ and $k$-handled tori. Formulas for Cauchy and Green's kernels are given in \cite{KraRyan2,BCKR}. 
A few years earlier, formulas for these kernels have been developed for spheres and hyperbolae as well, cf. \cite{lr,v}. In these works it has also been shown how many classical harmonic analysis techniques, such as for instance Calderon-Zygmund type operators (cf. \cite{lmq,m}) can be carried over from the classical Euclidean setting to these examples of spin manifolds. Actually, in the context of Riemannian spin manifolds, the Dirac operator arises as a special case of the Atiyah-Singer operator acting on sections of a spinor bundle, see also \cite{ER} for more details.  
\par\medskip\par
In \cite{KraRyan2} J. Ryan and the author also treated one particular example of a conformally flat manifold that is not orientable which hence is not spin, namely the real projective space $\R P^n$. The real projective arises by factoring out the unit sphere $S^n$ by the group $\mathbb{Z}_2$. 
\par\medskip\par
The aim of this paper is to carry over the previously developed techniques for a larger class of non-orientable conformally flat manifolds. Cylinders and tori, which arise by factoring out $\R^n$ by a discrete translation group, have natural non-orientable counterparts. First we look at a class of projective cylinders and tori, where we identify for example ``upper'' and ``lower'' parts with each other. This is translated into an additional symmetry structure on the manifold. 
\par\medskip\par
We explain how we can construct from the kernels for the cylinders and tori explicit kernels for this class of non-orientable manifolds. Similarly to the oriented case, we can consider distinct pin bundles in the non-oriented case. We show how arbitrary pin sections can be represented by these kernels. We also set up integral representation formulas for the sections in the kernel of the Dirac resp. Laplace operator on these manifolds. After that we discuss Hardy space decompositions in this context and shed some light on essential differences to the classical context. Then we set up integral formulas that permit us to compute the order of zeroes of sections on these manifolds.  
\par\medskip\par
Finally, we look at the M\"obius strip, the Klein bottle and some natural higher dimensional generalizations. These manifolds are also constructed by the same Kleinian group that generated the cylinders and tori, but with a different action that induces a switch of minus sign after one or several periods. This class of manifolds also play a crucial role in string and M-theory, for details see \cite{BBS}.  In turn the classical M\"obius strip in $\mathbb{R}^3$ finds its applications in mechanical engineering as well as in superconduction. 
\par\medskip\par    
In contrast to the class of non-orientable projective manifolds discussed previously, it is neither possible to construct monogenic sections  on the M\"obius strip nor on its higher dimensional generalizations that we discuss. This is an essential difference to the class of non-orientable manifolds discussed before. However, we are able to construct a harmonic Green's kernel that allows us to express harmonic sections on these manifolds in terms of boundary integrals over that kernel so that parts of the previously established machinery can be carried over to this context, too. 

\par\medskip\par
{\bf Acknowledgement}. The author is very thankful to Professor Dr. Karl Hofmann from the Technical University of Darmstadt (Arbeitsgruppe Algebra, Geometrie und Funktionalanalysis) for the fruitful discussions on the group theoretical construction methods of the Klein bottle. 

\section{Preliminaries}

For particular details about Clifford algebras see for instance \cite{p}. For basic Clifford analysis we refer the reader for example to  \cite{bds,GHS,GS2}. Let $\{e_1, e_2,\ldots, e_n\}$ be the standard basis of $\mathbb{R}^n$. ${\it Cl}_n$ denotes the associated real Clifford algebra in which $e_i e_j + e_j e_i = - 2 \delta_{ij}$ holds, where  $\delta_{ij}$ is the Kronecker symbol. 

Under this rule of multiplication each non-zero vector $x\in \R^{n}$ has a multiplicative inverse $x^{-1}=\frac{-x}{|x|^{2}}$. 
We also need the reversion anti-automorphism $\sim :Cl_{n}\rightarrow Cl_{n}:\sim e_{j_{1}}\ldots e_{j_{r}}=e_{j_{r}}\ldots e_{j_{1}}$. 

\par\medskip\par

The Euclidean Dirac operator is $D: = \sum\limits_{j=1}^n \frac{\partial }{\partial x_j} e_j$
and differentiable functions defined in open subsets of $\R^n$ with values in ${\it Cl}_n$ that are annihilated from the left (right) by the Dirac operator are called left (right) monogenic functions. 
The left and right fundamental solution to the $D$-operator is the Euclidean Cauchy kernel 
$G(x,y) = \frac{1}{\omega_{n}}\frac{x-y}{|x-y|^n}$ where $\omega_{n}$ is the surface area of the unit sphere $S^{n-1}$ lying in $\R^{n}$.
The Dirac operator factorizes the Laplacian, viz $D^2 = - \Delta$. For $n > 2$ the harmonic Green's kernel has the form $H(x,y) = \frac{1}{\omega_n(1-n)} \frac{1}{|x-y|^{n-2}}$. 

These kernels can be used to solve special boundary value problems. See for instance \cite{GS2,m} for more details.

\par\medskip\par

For all that follows it is crucial that the operators $D$ and $\Delta$ are invariant up to a conformal weight factor under all M\"{o}bius transformations.

In \cite{a} and elsewhere it is shown that any M\"{o}bius transformation $\psi(x)$ over $\R^{n}\cup\{\infty\}$ can be written in the form  $y=(ax+b)(cx+d)^{-1}$ where the coefficients $a,b,c,d$ are all products of vectors from $\R^{n}$, satisfying additionally
$a\tilde{c}$, $c\tilde{d}$, $d\tilde{b}$ and $b\tilde{a} \in \R^{n}$, and $a\tilde{d}-b\tilde{c} = \pm1$. 

Let $J_1(\psi,x) = \frac{\widetilde{cx+d}}{|cx+d|^{n}}$ and $J_2(\psi,x) = \frac{1}{|cx+d|^{n-2}}$. Here $|\cdot|$ represents the canonical extension of the Euclidean norm from $\R^n$ to $Cl_n$. 
If $f$ is a left monogenic (harmonic) function in the variable $y=\psi(x)=(ax+b) (cx+d)^{-1}$, then the function $J_{1}(\psi,x)f((ax+b) (cx+d)^{-1})$ or $J_{2}(\psi,x)f((ax+b) (cx+d)^{-1})$ is again left monogenic (or harmonic) now with respect to the variable $x$.  See \cite{r85,r2} and elsewhere. 

\par\medskip\par

\section{Spin and pin structures on conformally flat manifolds}

Conformally flat manifolds are in general $n$-dimensional manifolds that possess atlases whose transition functions are conformal maps in the sense of Gauss. For $n > 3$ the set of conformal maps coincide with the M\"obius transformations introduced previously. In the case $n=2$ the sense preserving conformal maps are exactly the holomorphic maps. So, under this viewpoint we can interpret conformally flat manifolds as higher dimensional generalizations of holomorphic Riemann surfaces. On the other hand, conformally flat manifold are precisely those Riemannian manifolds which have a vanishing Weyl tensor. 

\par\medskip\par

As mentioned for instance in the classical work by N. Kuiper \cite{Kuiper} concrete examples of conformally flat orbifolds can be constructed by factoring out a simple connected domain $X$ by a Kleinian group $\Gamma$ that acts discontinuously on $X$. In the cases  where $\Gamma$ is torsion free, the topological quotient $X/\Gamma$, consisting of the orbits of a pre-defined group action $\Gamma \times X \to X$, is endowed with a differentiable structure. We then deal with examples of conformally flat manifolds. 
    
\par\medskip\par

A classical way of obtaining pin or spin structures for a given Riemannian manifold is to look for a lifting of the principle bundle associated to the orthogonal group $O(n)$ or the special orthogonal group $SO(n)$, to a principle bundle for the pin group $Pin(n)$ or spin group $Spin(n)$. The $Spin(n)$ group is a subgroup of $Pin(n)$ of index $2$. $Spin(n)$ consists exactly of those matrices from $Pin(n)$ whose determinant equals $+1$. 
Now $Spin(n)$ is a double cover of $SO(n)$. So there is a surjective homomorphism $\theta: Spin(n)\rightarrow Spin(n)$ with kernel ${\mathbb{Z}}_{2}=\{\pm 1\}$. As explained in Appendix C of \cite{pp} this gives rise to a choice of two local liftings of the principle $SO(n)$ bundle to a principle $Spin(n)$ bundle. The number of different global liftings is given by the number of elements in the cohomology group $H^{1}(M,{\mathbb{Z}}_{2})$. See \cite{lm} and elsewhere for details. These choices of liftings give rise to different spinor bundles over $M$. Similarly, $Pin(n)$ is a double cover of $O(n)$ and analogously different pin bundles are obtained. For simplicity let us now first focus on the spin structures, because pin structures can be described analogously.

\par\medskip\par

In the case of a conformally flat manifold the fact that the group ${\mathbb{Z}}_{2}$ plays a central role in determining different spinor bundles can be seen immediately by noting that any M\"{o}bius transformation $y=\psi(x)$ can be either written as $(ax+b)(cx+d)^{-1}$ or as $(-ax-b)(-cx-d)^{-1}$. So in fact the conformal invariance of monogenic resp. harmonic functions that we previously mentioned is only correct up to a sign. Thus, the left monogenic function $f(y)$ is altered into $\pm J_1(\psi,x)f(\psi(x))$. This has an effect on constructing spinor bundles over a conformally flat manifold. Suppose that $M$ is a conformally flat manifold and $\mu_{2}\mu_{1}^{-1}:U_{1}\rightarrow U_{2}$ is a transition function arising from the atlas of $M$. So $U_{1}$ and $U_{2}$ are domains in $\R^{n}$ and $\psi=\mu_{2}\mu_{1}^{-1}$ is a M\"{o}bius transformation. Let us now consider the two trivial bundles $U_{1}\times Cl_{n}$ and $U_{2}\times Cl_{n}$. Given $u=\psi(x)\in U_{2}$ and $X\in Cl_{n}$ the pair $(u,X)\in U_{2}\times Cl_{n}$ may be identified with either $(x,J_1(\psi,x)X)$ or $(x,-J_1(\psi,x)X)$ in $U_{1}\times Cl_{n}$. If we can choose a suitable collection of signs on these local bundles which are globally compatible over $M$ then we have constructed a spinor bundle $E$ over $M$. In this case $M$ is called a conformally flat spin manifold. So in the same way that the choice of local liftings of principle $SO(n)$ bundles to local principle $Spin(n)$ bundles the $\mathbb{Z}_{2}$ dependence on the construction we have just given gives way to possible choices of spin bundles over conformally flat manifolds. 

Further it should be recalled, \cite{p}, that $Cl_{n}$ is the direct sum of several isomorphic minimal left ideals. These are often called spinor spaces. So in our construction of spinor bundles one can replace the Clifford algebra with one of these spinor spaces. 

Following \cite{ma,r85} we may now talk about monogenic (harmonic) sections.

\begin{definition}
Let $M$ be a conformally flat spin manifold with spinor bundle $E$. Then a section $f:M\rightarrow E$ is called a left monogenic (harmonic) section if locally $f$ reduces to a left monogenic (harmonic) function.
\end{definition}
The associated Laplacian on such a spin manifold is also called the spinorial Laplacian.
\par\medskip\par
Not all conformally flat manifolds are spin manifolds. A typical counterexample is the real projective space $\R P^n=S^n/\{\pm 1\}$. 

When $n$ is even, $\R P^{n}$ is no longer orientable. However, $\R P^n$ does still admit pin structures in these cases and hence pin  bundles similar to those described earlier in this section.

\section{Conformally flat cylinders and tori}

In this subsection we recall one basic construction method for the Cauchy and Green's kernel on the oriented conformally flat $n$-torus and conformally flat $k$-cylinders, which has been developed in \cite{KraRyan1,KraRyan2}. 
\par\medskip\par
To start, let us consider a general lattice  
$$
\Omega_n:= \mathbb{Z} v_1 + \cdots + \mathbb{Z} v_n
$$
where the elements $v_i$ $(i=1,...,n)$ are supposed to be $\mathbb{R}$-linearly independent vectors from $\mathbb{R}^n$. 
Next let us consider the group action 
$$
\Omega_n \times \mathbb{R}^n \to \mathbb{R}^n,\; v \circ x \mapsto x+v,\quad v \in \Omega,\;x \in \mathbb{R}^n.
$$
The invariance group of this lattice is the discrete translation group generated by the special M\"obius transformations $\psi(x) = x+v_i$ for $i=1,...,n$, so $J_1(\psi,x) \equiv \pm 1$ and $J_2(\psi,x) \equiv 1$. 
\par\medskip\par
Since $\Omega_n$ is a torsion free discrete Kleinian group, the topological quotient $T_n:=\mathbb{R}^n/\Omega_n$ consisting of the orbits of the above defined group action actually is a conformally flat manifold in $n$ real variables. It represents an $n$-dimensional torus. 
  \par\medskip\par
Of course any subgroup of a torsion free Kleinian group is again torsion free. Let $\Omega_k$ be a $k$-dimensional sublattice of $\Omega_n$. So, in particular for all $k=1,...,n-1$ the quotient sets   $C_k:=\mathbb{R}^n/\Omega_k$ which arise by the same group action are conformally flat manifolds, too. While $C_n=:T_n$ (the $n$-torus)  is compact, the other manifolds $C_k$ with $k=1,...,n-1$ are non-compact and have $n-k$ unbounded directions. We call them $k$-cylinders. The classical infinite cylinder is obtained by $\mathbb{R}^2/\mathbb{Z} e_1$. It arises from the group action $\mathbb{Z} \times \mathbb{R}^2 \to \mathbb{R}^2$ defined by $m \circ (x,y) \mapsto (x+m,y)$, where $m \in \mathbb{Z}$ and $(x,y) \in \mathbb{R}^2$. 
  
\par\medskip\par

Next, following the preceding paper \cite{KraRyan2}, the decomposition of the lattice $\Omega_k$ into the direct sum of the sublattices $\Omega_l:= \mathbb{Z} v_1 + \cdots + \mathbb{Z} v_l$ and $\Omega_{n-l}:= \mathbb{Z} v_{l+1} + \cdots + \mathbb{Z} v_k$ gives rise to $2^k$ distinct spinor bundles on $C_k$ --- that we will denote by $E^{(l)}$ --- by making the identification $(x,X) \leftrightarrow (x+\underline{m}+\underline{n}),(-1)^{m_1+\cdots+m_l}X)$ with $x \in \mathbb{R}^n,X \in Cl_n$. Here, we use the notation $\underline{m}= m_1 v_1 + \cdots+ m_l v_l \in \Omega_l$ and $\underline{n}:= n_{l+1} v_{l+1} + \cdots + n_k v_k \in \Omega_{k-l}$ where $m_1,\ldots,m_k \in \mathbb{Z}$.  In all that follows we denote the trivial bundle by $E^{(1)}$. 
Since the manifolds $C_k$ are orientable, we are dealing with examples of spin manifolds in this context here. 
\par\medskip\par
Notice that the different spin structures on a spin manifold $M$ are detected by the number of distinct homomorphisms from the fundamental group $\Pi_{1}(M)$ to the group ${\mathbb{Z}}_{2}$. In the case of the $k$-cylinder $C_k$ we have that $\Pi_{1}(C_k)={\mathbb{Z}}^{k}$. There are two homomorphisms of $\mathbb{Z}$ to $\mathbb{Z}_{2}$. The first is $\theta_{1}:{\mathbb{Z}}\rightarrow {\mathbb{Z}}_{2}:\theta_{1}(k) \equiv 0 \mod 2$ while the second is the homomorphism $\theta_{2}:{\mathbb{Z}}\rightarrow {\mathbb{Z}}_{2}:\theta_{2}(k) \equiv 1 \mod 2$. Consequently, there are $2^{k}$ distinct spin structures on $C_{k}$. Consequently the $n$-torus $C_{n}=T_{n}$ has $2^{n}$ distinct spin structures. $T_{n}$ is also an example of a Bieberbach manifold. Further details of spin structures on the $n$-torus and other Bieberbach manifolds can be found in \cite{f,mp,pf}.
\par\medskip\par
As we need the following construction in the next section, we briefly recall how we can construct monogenic and harmonic spinor sections, in particular the Cauchy and Green's kernel, on these manifolds. 

To proceed in this direction let us consider an open set $U \subset \mathbb{R}^n$.  
A function $f:U \to Cl_n$ that satisfies $f(x+\underline{m}+\underline{n}) = (-1)^{m_1+...+m_l} f(x)$ for all $v=\underline{m}+\underline{n} \in \Omega_l \oplus \Omega_{k-l}$ then naturally descends to a section on the $k$-dimensional cylinder $C_k(v_1,...,v_k) := \mathbb{R}^n/\Omega_k$ (including the $n$-torus when we put $k=n$) with values in the corresponding spinor bundle,  say $E^{(l)}$.  To do that one forms $f':=p_k(f)$, where $p_k: \mathbb{R}^n \to C_k$, $x \mapsto x':=x \mod C_k$ is the canonical projection from $\mathbb{R}^n$ down to the manifold $C_k(v_1,...,v_k)$.  

For the sake of simplicity we shall write $C_k$ instead of $C_k(v_1,...,v_k)$ when it is clear which basis vectors $v_1,...,v_k$ are considered. 

\par\medskip\par
 
The application of the projection map $p_k$ to the Dirac operator $D$ and to the Euclidean Laplacian $\Delta$ returns a first order differential operator $D'$ and a second order differential operator $\Delta'$ on the spin manifold $C_k$. The induced operators are the cylindrical (toroidal) Dirac and Laplace operators. 

\par\medskip\par
To construct the fundamental solution to these operators on these manifolds we basically sum the fundamental solution to the ordinary Euclidean Dirac or Laplace operator over the period lattice $\Omega_k$, taking into account the minus sign in correspondence to the spinor bundle that we consider. 

Recalling from \cite{KraRyan2}, for $k \le n-2$ the Cauchy kernel to the Dirac operator on $C_k$ with values in the spinor bundle $E^{(l)}$ can be expressed by 
$$
G^{(l)}_k(x',y')=p_k\left(\sum\limits_{\underline{m} \in \Omega_l} \sum\limits_{\underline{n} \in \Omega_{k-l}} (-1)^{m_1+\cdots+m_l}G(x-y+\underline{m}+\underline{n})\right).
$$
This series is majorized by the Eisenstein type series 
$$
\sum\limits_{\underline{\omega} \in \Omega_k} |\underline{\omega}|^{1-n}
$$ 
which is absolutely convergent for $k <n-1$. 
With a similar argument, for $k \le n-3$ one can express the Green's kernel to the Laplacian on $C_k$ with values in $E^{(l)}$ by 
$$
H^{(l)}_k(x',y')= p_k\left(\sum\limits_{\underline{m} \in \Omega_l} \sum\limits_{\underline{n} \in \Omega_{k-l}} (-1)^{m_1+\cdots+m_l} |x-y+\underline{m}+\underline{n})|^{2-n}\right).
$$
As described in \cite{KraRyan2} these sections indeed serve as the Cauchy or Green kernel in the corresponding Cauchy or Green's integral formula, respectively. They reproduce the values of interior points of monogenic (respectively harmonic)  spinor sections in a domain in terms of the given boundary data.  
\par\medskip\par
As proved in \cite{KraHabil,KraRyan2}, in the case $k=n-1$ the Cauchy kernel on $C_{n-1}$ can be expressed by the projection of the following modified monogenic  Eisenstein series 
$$
G^{(l)}_{n-1}(x,y) =   G(x-y)+\!\!\!\!\!\!\sum\limits_{(\underline{m},\underline{n}) \in \Omega_l \oplus \Omega_{n-1-l} \backslash\{(0,0)\} }\!\!\!\!\!\!\!\!\! (-1)^{m_1+\cdots+m_l}\Big[G(x-y+\underline{m}+\underline{n})-G(\underline{m}+\underline{n})\Big].
$$
The additional term $G(\underline{m}+\underline{n})$ guarantees that the series converges normally on $\R^n \backslash \Omega_k$. 
The same modification can be made for the harmonic kernel in the case $k=n-2$ to achieve the convergence of the series.  

\par\medskip\par

To round off we briefly recall from \cite{KraRyan2} how to construct a Cauchy kernel for the torus $T_n$. As we know for instance from Chapter 2.3 of \cite{KraHabil}, there does not exist a non-constant monogenic $n$-fold periodic function with only one point singularity of the order of the Cauchy kernel in one period cell. A non-constant $n$-fold periodic monogenic function must have either at least two point singularities of that order or a singularity of higher order in each period parallelepiped. Therefore, we only got a local Cauchy integral formula valid in sufficiently small domains. One can construct a non-constant and convergent monogenic $n$-fold periodic Eisenstein series whose projection then has two isolated point singularities at $a'$ and $b'$ of the order of the Cauchy kernel on the torus as follows. Take two points $a \not\equiv b \;{\rm mod}\; \Omega_n$ and consider   
\[ G(x-a) + G(x-b)+\sum\limits_{(\underline{m},\underline{n}) \in \Omega_l \oplus \Omega_{n-l} \backslash\{(0,0)\}}\Big[(-1)^{m_{1}+\ldots+m_{l}}(G(x-a+\underline{m}+\underline{n})\]
\[ +G(-a-\underline{m}-\underline{n})
+G(x-b+\underline{m}+\underline{n}) + G(-b-\underline{m}-\underline{n})\Big].\]
Due to this special coupling of terms in the series we obtain the further degree of convergence.  
In domains $D$ where either $a\not\in D$ or $b \not\in D$ and under some additional technical restrictions one can formulate a Cauchy integral formula in terms of such an Eisenstein series as kernel. For the particular technical details, we refer the reader to \cite{KraRyan2}. Similar modifications can be made for the harmonic case to deal with the cases $k=n-1$ and $k=n$.  

\par\medskip\par

In the following section we explain how we can adapt the construction method that we presented in this subsection in order to deal with some non-orientable counterparts of the manifolds considered here which can be constructed from the same Kleinian groups by which we constructed these cylinders and tori.

\section{Two classes of non-orientable conformally flat manifolds}

In this section we present the main results of our paper. 

\subsection{The geometric context}

To start, we recall that the oriented cylinder $C$ defined as the topological quotient $\mathbb{R}^2/\mathbb{Z}$ has several natural non-oriented counterparts. 
\par\medskip\par
First of all there is the projective cylinder which arises from the usual infinite cylinder by identifying the ``northern'' hemisphere with the ``southern'' hemisphere. In this case one has an additional symmetry relation on the original topological cylinder. This relation is mathematically expressed by $f(x,y)=f(x,-y)$ in addition to the usual periodic relation $f(x+m,y)=f(x,y)$ $(m \in \mathbb{Z})$. Recall that the standard cylinder is constructed by gluing the two vertices of a two-dimensional strip of length $1$ in a straight parallel way together without performing any twists. Here, one has besides the usual group action on the $x$-component of the form $\mathbb{Z} \times \mathbb{R}^2 \to \mathbb{R}^2$, $m\circ (x,y) = (x+m,y)$ an additional group action $\mathbb{Z}_2 \times \mathbb{R} \to \mathbb{R}$ on the second component, mapping $y$ and $-y$ on the same value.   

Functions on the projective standard cylinders consequently have an additional invariance under the action of the group $\{\pm 1 \} \cong \mathbb{Z}_2$ in the $y$-component.  
\par\medskip\par
A second classical non-oriented counterpart of the classical infinite cylinder $C$ is the M\"obius strip. The M\"obius strip has the property that a normal vector at a point $(x,y)$ is transformed into the opposite directed normal vector when returning to the same point after one period. Here we consider the same group $\Gamma=\mathbb{Z}$ and the same set $X=\mathbb{R}^n$, but one considers a different action of the form $\mathbb{Z} \times \mathbb{R}^2 \to \mathbb{R}^2$ where  $m\circ (x,y) = (x+m,(-1)^m y)$ ($m \in \mathbb{Z})$. 

In order to construct sections on the M\"obius strip one has to look at functions from the covering space $\mathbb{R}^2$ that satisfy  $f(x+m,y)=f(x,(-1)^m y)$ for $m \in \mathbb{Z}$. Only after two periods the normal vector at a surface point of the M\"obius strip again is transformed into itself in this geometry. Functions in $\R^2$ that descend to the standard M\"obius strip thus have period $2$ in the $x$-component.   
\par\medskip\par
Also the torus $T_2:=\mathbb{R}^2/\mathbb{Z}^2$ has such non-orientable counterparts. Its most famous candidate is the Klein bottle. As in the torus case, one considers $\Gamma = \mathbb{Z}^2$ and $X = \mathbb{R}^2$. However, the action is defined as 
$\mathbb{Z}^2 \times \R^2 \to \R^2:\;(m,n)\circ (x,y) \mapsto (x+m,(-1)^n y+n)$. A function on the Klein bottle thus descends from a function in $\mathbb{R}^2$ that satisfies the transformation rule $f(x+m,y+n)=f(x,(-1)^n y)$ for all $(m,n) \in \mathbb{Z}^2$. Here, we have a usual periodicity in the first component (like for the usual cylinder $C$). In contrast to the torus case we have a pseudo-periodicity with a minus sign switch depending on the parity of the second period in the second component. In contrast to the M\"obius strip, the minus sign switch here occurs in one of the periodic arguments.
\par\medskip\par
Notice that both the M\"obius strip and the Klein bottle can be obtained by gluing the same vertices of the fundamental domain of the associated one-dimensional  (respectively two-dimensional) translation group together that is used in the construction of the cylinder or  torus. However, one set of two opposite vertices of the fundamental domain is glued together with opposite orientation, causing a twist. This twist however destroys the orientability of the manifold. 
\par\medskip\par
In the $n$-dimensional setting we can construct similar non-oriented analogues of these manifolds from the oriented $k$-cylinders defined in the previous section by $C_k := \mathbb{R}^n/\Omega_k$ where $k \in \{1,...,n-1\}$. Here and in all that follows, $\Omega_k \subset \mathbb{R}^k$ denotes a $k$-dimensional lattice spanned entirely by $k$ $\mathbb{R}$-linearly independent reduced vectors $\underline{\omega}_1,...\underline{\omega}_k \in \mathbb{R}^k$. 

\par\medskip\par
{\bf Class A: Projective $k$-cylinders}
\par\medskip\par
Let $\underline{x}$ be a reduced vector from $\mathbb{R}^k$. Suppose that $\underline{\omega}:= m_1 \underline{\omega}_1 + \cdots + m_k \underline{\omega}_k$ is a vector from the lattice $\Omega_k$. Further let $p$ be an integer from the set $\{k+1,...,n\}$ in all that follows. 
\par\medskip\par
By identifying the tupel $$(\underline{x}+\underline{\omega},x_{k+1},..,x_p,..,x_n,X)$$ with $$(\underline{x},-x_{k+1},...,-x_p,...,x_{p+1},...,x_n,X)$$ we obtain an entire class of conformally flat manifolds denoted by ${\cal{M}}_{k,p}$. These manifolds are non-orientable and hence not spin manifolds anymore. 

\par\medskip\par
In the case not having any sign change at all, we again deal with the class of oriented conformally flat cylinders and tori which we briefly discussed in the previous section of this paper.  

\par\medskip\par

We can say more. Analogously, to the case of a spin manifold we can set up several distinct pin bundles. 
One pin bundle which is different from the trivial one given above is obtained by identifying the pair $$(\underline{x}+\underline{\omega},x_{k+1},...,x_n,X)\; {\rm with}\; (\underline{x}, -x_{k+1},...,- x_p,x_{p+1},...,x_n,-X).$$ Other distinct choices arise by again splitting the period lattice $\Omega_k$ into two sublattices $\Omega_k = \Omega_l \oplus \Omega_{k-l}$ where $1<l<k$. Writing an element $\underline{\omega} \in \Omega_k$ in the form $\underline{\omega} = \underline{m} + \underline{n}$ with $\underline{m} = m_1 \underline{\omega}_1+\cdots + m_l \underline{\omega}_l \in \Omega_k$ and $\underline{m} = m_{l+1} \underline{\omega}_{l+1}+\cdots + m_k \underline{\omega}_k \in \Omega_{k-l}$ gives rise to again consider the following identifications 
$$(\underline{x}+\underline{\omega},x_{k+1},...,x_n,X) \; {\rm with}\; (\underline{x}, -x_{k+1},...,-x_p,x_{p+1}..., x_n,(-1)^{m_1+\cdots+m_l}X).$$   

\par\medskip\par

{\bf Class B: Higher dimensional M\"obius strips and generalizations of the Klein bottle}

\par\medskip\par

Similarly to the classical case in three dimensions we can introduce higher dimensional analogues of the M\"obius strip as the topological quotient consisting of the orbits of the group action 
$$
\Omega_k \times \mathbb{R}^n \to \R^n
$$
where the action is defined by 
$$
\underline{\omega} \circ x := (\underline{x}+\underline{\omega},x_{k+1},...,x_{n-1},{\rm sgn}(\underline{\omega}) x_n).
$$
Here, for $\underline{\omega} = m_1 \underline{\omega}_1 + \cdots + m_k\underline{\omega}_k$ we write    
${\rm sgn}(\underline{\omega}) = \left\{ \begin{array}{cc} 1 & {\rm if}\; \underline{\omega} \in 2 \Omega_k \\ -1 & {\rm if}\; \underline{\omega} \in \Omega_k \backslash 2 \Omega_k \end{array} \right.$. 
To be more explicit this action describes manifolds of the form  
$$
{\cal{M}}_k^{-} = \mathbb{R}^n/\sim
$$
where $\sim$ is now defined by the map $$(\underline{x}+\underline{\omega},x_{k+1},...,x_{n-1},x_n) \mapsto (\underline{x},x_{k+1},...,x_{n-1},{\rm sgn}(\underline{\omega})x_n).$$
In the case $n=2,k=1$ we re-obtain the classical two-dimensional M\"obius strip. Notice that - again due to the switch of the minus sign - the manifolds ${\cal{M}}_k^{-}$ are not spin manifolds.  As in the previous cases dealt in Class A we can of course again involve sign changes in several of the remaining directions $k+1,...,n$. Alternatively, we can re-define the expression ${\rm sgn}$ in such a way that ${\rm sgn}(\underline{\omega}) = 1$ if $\sum_{i=1}^k m_i \in 2 \mathbb{Z}$ and ${\rm sgn}(\underline{\omega}) = -1$ otherwise. 
\par\medskip\par
When we involve minus signs in one or several components on which the periodicity lattice acts (that means on components $j$ with $j \le k$) we obtain natural higher dimensional generalizations of the Klein bottle. To leave it simple we shall consider a $k$-dimensional normalized lattice of the form $\Omega_k = \Omega_{k-1} + \mathbb{Z} e_k$ with $\Omega_{k-1} \subset \R^{k-1}$ and 
a minus sign switch in the last periodicity component $k$ only, which only depends on the last period. Notice that any arbitrary $k$-dimensional lattice from $\R^k$ can be transformed into a lattice of this form. This can be achieved by applying one elementary rotation and one dilation. 

These cases are described by group actions of the form   
$$
\Omega_k \times \mathbb{R}^n \to \R^n
$$
where the action is defined by 
$$
\underline{\omega} \circ x := (\underline{x}+\sum_{i=1}^{k-1}m_i \underline{\omega}_i+(-1)^{m_k} m_k e_k,x_{k+1},...,x_{n-1}, x_n)
$$
The quotient manifolds can then be described in terms of 
$$
{\cal{K}}_k := \mathbb{R}^n/\sim^*
$$
where 
$\sim^*$ is now defined by the map $$(\underline{x}+\sum_{i=1}^{k-1}\underline{\omega}_i+m_k e_k,x_{k+1},...,x_{n-1},x_n) \mapsto (x_1,...,x_{k-1},(-1)^{m_k}x_k,x_{k+1},...,x_{n-1},x_n).$$
In the case where $k=n=2$ we re-obtain the classical Klein bottle. 

Notice that the manifolds ${\cal{K}}_k^{-}$ can also be constructed by gluing finitely many conformally flat manifolds together. This is another argument for ${\cal{K}}_k^{-}$ being conformally flat, cf. \cite{Ry2003}.   

More generally, we may consider minus sign changes in several periodicity components. In addition to that we can consider minus sign switches both in some of the periodicity components $j \in \{1,...,k\}$ as well as in some of the remaining components. The associated manifolds then represent a complicated mixture of the geometries of the M\"obius strip and of the Klein bottle.     
\par\medskip\par
  
Notice that in all these different contexts one again can set up several distinct pin bundles (different from the trivial one constructed above) by making the different identifications with the minus signs as in the projective cases. Again, this can be done by considering different decompositions of the period lattice $\Omega_k$. As the construction principle is clear after our explanations in the preceding subsection we leave a detailed description of the other pin bundles to the reader here.

\subsection{Construction of the Cauchy and Green's kernel, corresponding integral formulae and applications}

Let us first treat the manifolds from Class $A$. 

\par\medskip\par

For simplicity let us look in the following at the trivial pin bundle of the manifold ${\cal{M}}_{k,p}$  where we simply identify 
$$
(\underline{x}+\underline{\omega},x_{k+1},...,x_p,x_{p+1},...,x_n,X)$$ with $$(\underline{x}, -x_{k+1},...,- x_p,x_{p+1},...,x_n,X).
$$ 
Again we denote by $G_k(x',y')$ the monogenic Cauchy kernel for the Dirac operator and by $H_k(x',y')$ the harmonic Green's kernel for the Laplace operator on the infinite $k$-cylinder $C_k$ associated with the trivial bundle. To leave it simple we assume that in the monogenic case we have $k < n$. In the harmonic setting we assume that $k < n-1$.   
\par\medskip\par
Now we can express the Cauchy kernel of the Dirac operator and the Green's kernel of the Laplace operator on the manifold ${\cal{M}}_{k,p}$ by a finite superposition of the fundamental sections $G_k$ or $H_k$, respectively, taking care of the additional minus signs in the coordinate directions $k+1,...,p$. More concretely, we may establish that 
\begin{theorem}
The Cauchy kernel of the Dirac operator on ${\cal{M}}_{k,p}$ associated to the trivial pin bundle can be expressed in the form 
\begin{eqnarray*}
G_{{\cal{M}}_{k,p}}(x'',y'')&=&P_p \Bigg(  \sum\limits_{\varepsilon_{k+1},...,\varepsilon_{p} \in \{\pm 1\}} \!\!\!\!\!\!\!\! G_k(\underline{x'}-\underline{y'},\varepsilon_{k+1} (x'_{k+1}-y'_{k+1}),...\\
& &\quad\quad\quad\quad\quad ..., \varepsilon_{p} (x'_p-y'_p),(x'_{p+1}-y'_{p+1}),...(x'_n-y'_n))\Bigg). 
\end{eqnarray*}
Here, $P_p$ denotes the projection from the fully infinite cylinder $C_k$ down to the projective cylinder ${\cal{M}}_{k,p}$. 
In the cases of the other pin bundles mentioned before, we need to add the corresponding minus sign in the sum in front of the multiperiodic expression $G_k$.   
\end{theorem}
To prove this statement we first recall that in the preceding works \cite{KraHabil,KraRyan1} it has already been shown that the subseries $G_k(x,y)$ are well-defined and normally convergent in $\mathbb{R}^n \backslash \Omega_k$. Furthermore, they are invariant under the lattice transformations $G_k(x+v_i,y) = G_k(x,y+v_i)=G_k(x,y)$ for all $v_i \in \Omega_k$, and therefore for all $v \in \Omega_k$. The summation over all $\varepsilon_{k+1},...,\varepsilon_{p} \in \{\pm 1\}$ is a finite sum of expressions $G_k(x,y)$. Consequently, the total series again converges normally in $\mathbb{R}^n \backslash \Omega_k$ and is $\Omega_k$-invariant.  Due to the extension of the summation over all possible sign combinations in the coordinates $k+1,...,p$  the outcome then is also invariant under transformations of the form 
$f(\underline{x}-\underline{y},(x_{k+1}-y_{k+1}),...,(x_p-y_p),...,(x_n-y_n))$ $\mapsto
f(\underline{x}-\underline{y},-(x_{k+1}-y_{k+1}),..., -(x_p-y_p),(x_{p+1}-y_{p+1}),...,(x_n-y_n)),$ 
too. Therefore, the projection $P_{k,p}:=P_p \circ p_k$ of the entire series descends to a well-defined left monogenic section on ${\cal{M}}_{k,p}$. On the manifold ${\cal{M}}_{k,p}$ the kernel then has exactly one point singularity of the order of the usual Cauchy kernel. As a sum of left monogenic functions, it is again left monogenic except at this point singularity. 
As a consequence of the usual Cauchy theorem valid in the universal covering space $\R^n$, cf. \cite{bds}, which states that the oriented  boundary integral over a left monogenic function vanishes identically, we can deduce a properly adapted version of Cauchy's integral formula on the manifold from which then follows that $G_{{\cal{M}}_{k,p}}(x'',y'')$ actually is the Cauchy kernel. 

\par\medskip\par

To get there, let us first consider a strongly Lipschitz hypersurface $S$ that lies completely in that part of the standard fundamental period parallelepiped of the lattice $\Omega_k$, where all $x_{k+1},...,x_p > 0$, i.e. inside of $D:=\{x \in \mathbb{R}^n \mid 0 < x_i < m_i \;\;(i=1,...k)\;{\rm and}\; x_j > 0\;\;(k+1 \le j \le p)\}$. 
Strongly Lipschitz means that locally the hypersurface is the graph of a Lipschitz function, and that globally the local Lipschitz constants are bounded. Suppose also that $V$ is a domain lying in $D$ and that $S$ bounds a subdomain $W$ of $V$ such that $W \cup S \subset V$. Now take an element $y \in W$. As a consequence of the usual Cauchy formula and Cauchy's theorem in $\mathbb{R}^n$, we obtain that if $f:V \to Cl_n$ a left monogenic function, then 
\begin{eqnarray*}
f(y) &=& \! \int\limits_S \! \Big(\sum\limits_{\varepsilon_{k+1},...,\varepsilon_{p} \in \{\pm 1\}}\sum\limits_{\underline{\omega} \in \Omega_k} G(\underline{x}-\underline{y}+\underline{\omega},\varepsilon_{k+1} (x_{k+1}-y_{k+1}),... \\ & & \quad\quad\quad...,\varepsilon_p (x_p-y_p),(x_{p+1}-y_{p+1}),...,(x_n-y_n))\Big) \times \\ & & \quad\quad\quad \times \;\; n(x) f(x) d\sigma(x),
\end{eqnarray*}
where $n(x)$ stands for the unit exterior normal vector to $S$ at $x$ and $\sigma$ denotes the usual Lebesgue measure on $S$.  
\par\medskip\par
As $S,V,W$ entirely lie in $D$, the projection map $p_k$ induces a Lipschitz surface $S'$, a domain $V'$ and a subdomain $W'$ of $V'$ that lie entirely on that part of the infinite cylinder $C_k$ where additionally $x_{k+1},...,x_p > 0$. So, the original sets $V,W$ descend to well-defined domains $V'',W''$ and the surface $S$ to a well-defined strongly Lipschitz hypersurface $S''$ on the projective cylinder ${\cal{M}}_{k,p}$. 
As a consequence of the usual Cauchy's theorem from \cite{bds} we can now formulate 
\begin{theorem}\label{cauchy}(Cauchy integral formula).\\
Suppose that $S''$ is strongly Lipschitz hyperface bounding a subdomain $W''$ of a domain $V''$ lying on the manifold ${\cal{M}}_{k,p}$. Let $f'':V'' \to E^{1}$ be a left monogenic section on ${\cal{M}}_{k,p}$ with values in the trivial pin bundle $E^{1}$. Suppose that $y'' \in W''$. Then 
$$
f''(y'') = \int\limits_{S''} G_{{\cal{M}}_{k,p}}(x'',y'') dP_{k,p}(n(x)) f''(x'') d\sigma''(x''),
$$  
where $dP_{k,p}$ stands for the derivative of the projection map $P_{k,p} = P_p \circ p_k$. 
\end{theorem} 
\begin{remark} We can replace $S''$ by a nullhomologous $(n-2)$-dimensional cycle $\Gamma''$ parameterizing an
$(n-2)$-dimensional surface of an $n-1$-dimensional simply connected
domain inside of $G''\subset {\cal{M}}_{k,p}$. In this more general setting Cauchy's integral formula then takes the form 
$$
w_{\Gamma''}(y'')f''(y'') = \int\limits_{S''} G_{{\cal{M}}_{k,p}}(x'',y'') dP_{k,p}(n(x)) f''(x'') d\sigma''(x''),
$$   
where $w_{\Gamma''}(y'')$ is the wrapping number counting how often the cycle $\Gamma''$ wraps around the point $y''$. To derive the generalized homological Cauchy integral formula from the simple Cauchy integral formula one needs to apply the same arguments as used in the Euclidean case. This is described in detail in {\rm \cite{GHS}, Chapter 12.3.4}.
\end{remark}
In view of the well-known Almansi-Fischer decomposition theorem, cf. \cite{bds}, we can now directly set up a Green's integral formula for harmonic sections on these manifolds. To proceed in this direction, we need to construct the Green's kernel for harmonic sections on ${\cal{M}}_{k,p}$. We can obtain the Green's kernel on ${\cal{M}}_{k,p}$ as a finite sum of harmonic Green's kernels that we constructed for the fully infinite cylinders $C_k$. By similar arguments as given above the harmonic Green's kernel for harmonic sections on ${\cal{M}}_{k,p}$ with values in the trivial pin bundle $E^{1}$ can be expressed by 
\begin{eqnarray*}
H_{{\cal{M}}_{k,p}}(x'',y'')&=&P_p \Bigg(  \sum\limits_{\varepsilon_{k+1},...,\varepsilon_{p} \in \{\pm 1\}} \!\!\!\!\!\!\!\! H_k(\underline{x'}-\underline{y'},\varepsilon_{k+1} (x'_{k+1}-y'_{k+1}),\\
& &\quad\quad\quad\quad\quad ..., \varepsilon_{p} (x_p-y_p),(x'_{p+1}-y'_{p+1}),...(x'_n-y'_n))\Bigg) 
\end{eqnarray*}
where $H_k(x',y')$ is the harmonic Green's kernel for harmonic sections on $C_k$ with values in the trivial pin bundle. This tool in hand, we can establish, after having applied the Almansi-Fischer decomposition, the following 
\begin{theorem} (Green's integral formula).\\
Let $k < n-1$. Suppose that $V'',W'',S''$ and $y''$ are as in Theorem~\ref{cauchy}. Now let $f'':V''\to E^{1}$ be a harmonic section on ${\cal{M}}_{k,p}$. Then 
\begin{eqnarray*}
f''(y'') &=& \int\limits_{S''} G_{{\cal{M}}_{k,p}}(x'',y'') dP_{k,p}(n(x)) f''(x'') d\sigma''(x'') \\
         &+& \int\limits_{S''} H_{{\cal{M}}_{k,p}}(x'',y'') dP_{k,p}(n(x)) D''[f''(x'')] d\sigma''(x'')
\end{eqnarray*}
where $D''$ denotes the Dirac operator on ${\cal{M}}_{k,p}$ induced by the projection $P_{k,p}$. 
\end{theorem} 
Here, we see that the knowledge of the Cauchy kernel for the Dirac operator is indeed useful to study harmonic sections on these manifolds. In the cases $k=n-1,k=n$ we can also set up a local version of the Green's integral formula, which however then involves a Cauchy kernel with more than one singularity on the whole manifold ${\cal{M}}_{k,p}$, similarly to that on the $n-1$-cylinders and $n$-tori discussed in \cite{KraRyan2,KraRyan3}. 
\par\medskip\par
We proceed to give some further applications.  
\par\medskip\par
{\bf Hardy space decompositions}. Let $q$ be a real parameter within the range $(1,\infty)$. Let us suppose that $S$ is a strongly Lipschitz hypersurface that bounds a subdomain $W$ of $V$ such that $W \cup S \subset V$. Let us first additionally assume that $V \subset D$. In this case it is straightforward to deduce from our previous calculations for the fully infinite cylinders $C_k$, presented in \cite{KraRyan1}, that one obtains the usual Hardy space decomposition for $L^q(S'')$ in the form 
$$
L^q(S'') = H^q({S''}^+) \oplus H^q({S}''^-).
$$
Here ${S''}^+=P_{k,p}(W)$, ${S''}^-={\cal{M}}_{k,p} \backslash (P_{k,p}(W) \cup S'')$ and $H^q(S{''}^{\pm})$ denotes the Hardy $q$-space of left monogenic sections defined on ${S''}^{\pm}$ with non-tangential limits on $S''$ with values in $L^q(S'')$. 
\par\medskip\par
However, we shall now assume that the domain $V$ is such that for each $$x:=(x_1,...,x_k,x_{k+1},...,x_p,...,x_n) \in V$$ the element 
$$\tilde{x}:=(x_1,...,x_k,-x_{k+1},...,-x_p,x_{p+1},...,x_n) \in V,$$ too.
If we now assume that the hypersurface $S$ is such that  $\tilde{S} = S$ where again 
$$
\tilde{S} := \{s \mid \tilde{s} \in S\},
$$
then both $y$ and $\tilde{y}$ belong to $V$. This fact implies that after applying the projection $P_{k,p}$ we obtain on the manifold ${\cal{M}}_{k,p}$ that 
$$
\int\limits_{S''} G_{{\cal{M}}_{k,p}}(x'',y'') dP_{k,p}(n(x)) f''(x'') d\sigma''(x'') = 2 f''(y'').
$$ 
Returning to the covering space $\mathbb{R}^n$, let us now consider a function $\eta:S \to Cl_n$ which belongs to $L^q(S,Cl_n)$ with $1 < q < \infty$ satisfying $$
f(\underline{x}+\underline{\omega},x_{k+1},...,x_n)= f(x_1,...,x_k,x_{k+1},...,x_n) \quad \quad \forall \underline{\omega} \in \Omega_k
$$ 
and 
$$
f(\underline{x},-x_{k+1},...-x_p,x_{p+1},...,x_n)=f(\underline{x},x_{k+1},...,x_n).
$$
Let us consider a piecewise $C^1$ path $y(t) \in V$ which approaches $w \in S$ non-tangentially as $f$ tends to $0$. Then 
\begin{eqnarray*}
& & \lim\limits_{t \to 0} \int\limits_S \! \Big(\sum\limits_{\varepsilon_{k+1},...,\varepsilon_{p} \in \{\pm 1\}}\sum\limits_{\underline{\omega} \in \Omega_k} G(\underline{x}-\underline{y}(t)+\underline{\omega},\varepsilon_{k+1} (x_{k+1}-y_{k+1}(t)),... \\ & & \quad\quad\quad ....,\varepsilon_p (x_p-y_p(t)),(x_{p+1}-y_{p+1}(t)),...,(x_n-y_n(t)))\Big) \times \\ & & \quad\quad\quad \times \;\; n(x) \eta(x) d\sigma(x)
\end{eqnarray*}
evaluates almost everywhere to 
\begin{eqnarray*}
& & \frac{1}{2} \eta(w) + P.V. \int\limits_S \! \Big(\sum\limits_{\varepsilon_{k+1},...,\varepsilon_{p} \in \{\pm 1\}}\sum\limits_{\underline{\omega} \in \Omega_k} G(\underline{x}-\underline{z}+\underline{\omega},\varepsilon_{k+1} (x_{k+1}-z_{k+1}),... \\ & & \quad\quad\quad ....,\varepsilon_p (x_p-z_p),(x_{p+1}-z_{p+1}),...,(x_n-z_n))\Big) \times \\ & & \quad\quad\quad \times \;\; n(x) \eta(x) d\sigma(x)
\end{eqnarray*}
When we turn to the projective cylinder ${\cal{M}}_{k,p}$ we are now obliged to consider two paths $y(t)$ and $\tilde{y}(t)$ in the covering space $\mathbb{R}^n$. Consequently the limit
$$
\lim\limits_{t \to 0} \int\limits_{S''} G_{{\cal{M}}_{k,p}}(x'',y(t)'') dP_{k,p}(n(x)) {\eta}''(x'') d\sigma''(x'')
$$ 
is equal to the value 
$$
2 P.V. \int\limits_{S''} G_{{\cal{M}}_{k,p}}(x'',z'') dP_{k,p}(n(x)) {\eta}''(x'') d\sigma''(x'')
$$  
where $z'':=P_{k,p}(z)$ as well as $\eta'':=P_{k,p}(\eta)$. Here, we observe that the usual Hardy space decomposition as mentioned above does not always occur. It does occur, when $S$ lies in a positivity domain of the type $D$.

\begin{remark} In {\rm \cite{KraRyan1}} we also introduced convolution operators of Calderon-Zygmund type acting on $L^q$ spaces of special hypersurfaces on the fully infinite oriented cylinders $C_k$.  Furthermore, we introduced proper analogues of operators of LMS type and Poisson-, Szeg\"o-, Bergman- and Kerzman-Stein kernels. All these operators directly carry over to the context of the manifolds considered here, by adapting the kernel in the way that we presented here.
\end{remark}
{\bf Applications to order theory}. A further consequence of Cauchy's integral formula is the following version of an order formula for isolated zeroes of left monogenic sections on this class of manifolds. 

To proceed in this direction we first need to introduce the notion of the order of an isolated zero of a left monogenic section on the projective cylinder manifold ${\cal{M}}_{k,p}$. Again, for simplicity, we restrict ourselves to treat the trivial pin bundle, since the others can be treated analogously by making the proper adaptation in the series representation of the Cauchy kernel.  
\par\medskip\par
Let $c \in D$ and suppose that $\varepsilon > 0$ is sufficiently small so that the open ball $B(c,\varepsilon)$ lies completely in $D$. Next denote by $B''(c'',\varepsilon)$ the projection of that ball to the manifold  ${\cal{M}}_{k,p}$, i.e. $B''(c'',\varepsilon) := P_{k,p}(B(c,\varepsilon))$. Now we may introduce 
\begin{definition}
\label{ord}
Suppose that $G'' \subseteq {\cal{M}}_{k,p}$ is an 
open set and that $g'': G'' \rightarrow E^{(1)}$ is a left monogenic section on ${\cal{M}}_{k,p}$ with values in the trivial pin bundle $E^{(1)}$. Next assume that $c'' \in G''$ is an isolated zero on ${\cal{M}}_{k,p}$ so that there is a $\delta > 0$ (with $\delta < \varepsilon)$ such that the closure of the projection of the ball 
$\overline{B''}(c'',\delta) \subset G''$ and $g''|_{\overline{B''}(c'',\delta)
\backslash\{c''\}} \neq 0$. Then the integer  
$$
{\rm ord}(g'';c'') := \int\limits_{g''(\partial B''(c'',\delta))} G_{{\cal{M}}_{k,p}}(x'',0'') dP_{k,p}(n(x)) d\sigma''(x'')
$$
is called the order of $g''$ at $c''$.  
\end{definition}
\begin{remark}
This way to introduce the order of a zero is compatible with the classical notion of the topological mapping degree given in {\rm \cite{AH}} for harmonic functions in $\R^n$.  
\end{remark}
In order to show that ord $(g'';c'')$ really is an integer, we now need to apply the Cauchy 
integral formula that we now have established for the manifold ${\cal{M}}_{k,p}$. Let us put in Theorem~\ref{cauchy} 
$f''\equiv 1$ and replace $y''$ by $g''(c'')$ and $\partial S''$
by $g''(\partial B''(c'',\delta))$. This leads to
$$
\int\limits_{g''(\partial B''(c'',\delta))} 
G_{{\cal{M}}_{k,p}}(x'',\underbrace{g''(c'')}_{=0})dP_{k,p}(n(x)) d\sigma''(x'')=w_{g''(\partial B''(c'',\delta))}(\underbrace{g''(c'')}_{=0}).
$$
The value 
ord$(g'',c'')=w_{g''(\partial B''(c'',\delta))}(0)$ represents the integer
counting how often the image of $B''$ under $g''$ around the isolated zero
wraps around zero.

Similarly to the Euclidean case, one can replace the projection of the ball in Definition~\ref{ord} again by a nullhomologous $(n-2)$-dimensional cycle parameterizing an 
$(n-2)$-dimensional surface of an $n-1$-dimensional simply connected
domain inside of $G''\subset {\cal{M}}_{k,p}$ which contains the isolated zero $c''$ in its interior and no further zeroes neither in its interior nor on its boundary.

Let us now assume that $G''$ is a domain. 
Since $g'':G'' \rightarrow E^{(1)}$ is continuously differentiable (because of the left monogenicity), we
can then again apply the general transformation rule for differential forms, cf. for example \cite{KraHabil}).   
For the oriented surface differential on the manifold ${\cal{M}}_{k,p}$, abbreviated by 
$d\Sigma''(x''):=dP_{k,p}(n(x)) d\sigma''(x'')$ we obtain the following formula   
\begin{equation}
\label{zoell1cyl} 
d\Sigma''(g''(x''))=[(Jg'')^*(x'')] * [d\Sigma''(x'')].
\end{equation}
Here, $Jg''$ denotes the Jacobi matrix of $g''(x'')$  and $(Jg'')^*$ its adjoint. The symbol $*$ between the brackets indicates again the matrix-vector multiplication.  
This allows us to rewrite the expression of the order in the following form.  
\begin{theorem}
\label{argumenttorus}
Let $G''\subseteq {\cal{M}}_{k,p}$ be a domain. Let $g'':G''\rightarrow E^{1}$ be a  
left monogenic section in $G''$ and suppose that $c''\in G''$  is an isolated zero. Then, under the same
conditions for $B''(c'',\delta)\subset G''$ as in Definition~\ref{ord},
\begin{equation}
\label{argcylinder1}
{\rm ord}(g'';c'')=
\int_{\partial B''(c'',\delta)}
G_{{\cal{M}}_{k,p}}(x'',0'')
(g''(x'')) [(Jg'')^{*}(x'')]*[d\Sigma''(x'')].
\end{equation}
\end{theorem}  
\begin{remark} By means of formula~(\ref{argcylinder1}) one can then establish in a similar way as in the Euclidean case treated in Chapter~1.5 of {\rm \cite{KraHabil}} an explicit argument principle as well as a Rouche's theorem adapted to the proper metric on ${\cal{M}}_{k,p}$.  In this context notice that whenever $\Gamma''$ is an  $(n-2)$-dimensional cycle on ${\cal{M}}_{k,p}$ which has atmost a countable number of points $z''$ at which $\det Jg''(x'') = 0$, then we may rewrite the expression $[(Jg'')(x'')]*[d\Sigma(x'')]$ in the form $\det Jg''(x'') [((Jg'')^{-1})^\top(x'')]*[d\Sigma''(x'')]$.
\end{remark}

\begin{remark} In the extreme case where one puts $k=0$ in which we identify the pair $(x_1,...,x_n,X)$ with $(-x_1,\ldots, -x_p,x_{p+1},...,x_n,X)$ we are dealing with real projective manifolds without periodicity conditions. The latter shall be denoted by ${\cal{M}}_{0,p}$. In this limit case we simply have to replace the function  
$$G_k(\underline{x}-\underline{y},\varepsilon_{k+1} (x_{k+1}-y_{k+1}),..., \varepsilon_{p} (x_p-y_p),(x_{p+1}-y_{p+1}),...,(x_n-y_n))$$ 
by the function $$G(\varepsilon_{1} (x_{1}-y_1),..., \varepsilon_{p} (x_p-y_p),(x_{p+1}-y_{p+1}),...,(x_n-y_n))$$ and the Cauchy kernel to the Dirac operator on ${\cal{M}}_{0,p}$ is simply given by 
$$
\sum\limits_{\varepsilon_{1},...,\varepsilon_{p} \in \{\pm 1\}}  G(\varepsilon_{1} (x_{1}-y_{1}),..., \varepsilon_{p} (x_p-y_p),(x_{p+1}-y_{p+1}),...,(x_n-y_n)).
$$
Notice that this is a finite series only. Similarly, one constructs the corresponding harmonic kernel for the Laplacian. 
\end{remark}

Now we turn to the manifolds of Class $B$. Here, we are confronted with essential differences. Differentiable functions on the manifolds ${\cal{M}}_{k}^-$ arise from differentiable functions defined in $\mathbb{R}^n$ that satisfy $$f(\underline{x}+\underline{\omega},x_{k+1},...,x_{n-1},x_n)=f(\underline{x},x_{k+1},...,x_{n-1},{\rm sgn}(\underline{\omega})x_n)$$ for all $\underline{\omega} \in \Omega_k$. However, if a non-constant function $f(x_1,...,x_{n-1},x_n)$ is annihilated by the Euclidean Dirac operator $D$ associated with the vector variable $x=(x_1,...,x_n)$, then the function $$g(x_1,...,x_{n-1},x_n):=f(x_1,...,x_{n-1},-x_n)$$ is not anymore in the kernel of operator $D$ where differentiation is meant with respect to the same variables $x_1,...,x_n$. As a consequence one directly obtains
\begin{theorem}
The only left monogenic functions on the manifolds ${\cal{M}}_k^-$ are constants. 
\end{theorem} 
However, we can construct harmonic functions, or more generally, harmonic sections on these manifolds. Differentiating twice with respect to each variable $x_1,...,x_n$, it follows that if a function $f(x_1,...,x_n)$ is harmonic with respect to the variables $x_1,...,x_n$, then so is $f(x_1,...,-x_n)$ (with respect to the same variables $x_1,...,x_n$). 

To obtain the Green's kernel for harmonic sections on these manifolds, we can adapt the kernel formula that we obtained for the $k$-cylinders, namely by inserting a minus sign in the $n$-th coordinate of the expressions appearing in the sum that are associated to those lattice points with ${\rm sgn}(\underline{\omega}) = -1$. More precisely, for the trivial pin bundle we obtain 
\begin{theorem} Let $k < n-2$ and let ${\cal{M}}_k^-$ be the manifold constructed above.  Then the harmonic Green's kernel for the trivial pin bundle is given by
$$
H^{-}_k(x',y') =  p_-\Bigg(\sum\limits_{\underline{\omega} \in \Omega_k}  |(\underline{x}-\underline{y}+\underline{\omega})+(x_{k+1}-y_{k+1}) e_{k+1} + ... + (x_{n-1}-y_{n-1}) e_{n-1} 
$$
$$
 + {\rm sgn}(\underline{\omega}) (x_n-y_n) e_n |^{2-n}\Bigg)
$$
where $p_-$ again stands for the projection that maps a point $x \in \mathbb{R}^n$ to $x'=x \mod \sim$ belonging to the manifold ${\cal{M}}_k^-$.  In the case $k=n-2$, the sum is substituted by the expression 
$$
 |x-y|^{2-n}+\sum\limits_{\underline{\omega} \in \Omega_k \backslash \{\underline{0}\}}  \Big(|(\underline{x}-\underline{y}+\underline{\omega})+(x_{k+1}-y_{k+1}) e_{k+1} + ... + (x_{n-1}-y_{n-1}) e_{n-1}$$
$$
+ {\rm sgn}(\underline{\omega}) (x_n-y_n) e_n |^{2-n} - |\underline{\omega}|^{2-n}\Big).
$$
\end{theorem} 
To the proof one observes that in the cases $k < n-2$ the series 
$$f(x):= \left(\sum\limits_{\underline{\omega} \in \Omega_k}  |(\underline{x}+\underline{\omega})+x_{k+1} e_{k+1} + ... + x_{n-1} e_{n-1} + {\rm sgn}(\underline{\omega}) x_n e_n |^{2-n}\right)$$
indeed satisfies 
$$
f(\underline{x}+\underline{\eta},...,x_{n-1},x_n)=f(\underline{x},...,x_{n-1},{\rm sgn}(\underline{\eta}) x_n)
$$
for all $\eta \in \Omega_k$ which follows by a direct rearrangement argument of the series. The same holds for $k=n-2$, but the rearrangement argument is a bit more sophisticated. (See \cite{KraHabil} where the oriented cylinder case is treated.) 

So, its projection defines a well-defined harmonic pin section on ${\cal{M}}_k^-$. On ${\cal{M}}_k^-$ it has precisely one isolated point singularity of the order of the usual harmonic Green's kernel. The argumentation that this kernel actually reproduces the harmonic sections on ${\cal{M}}_k^-$ can be done along the same lines of the argumentation that we gave for the manifolds of Class A in the preceding part of this subsection. The main ingredient again is the fact that the boundary integral over a function vanishes when it is harmonic in all its interior points. 

\par\medskip\par 
 
To round off, we finally turn to the manifolds ${\cal{K}}_k^-$ that generalize the Klein bottle to higher dimensions. As in the case of the higher dimensional M\"obius strips there are no non-constant left monogenic sections on ${\cal{K}}_k^-$ for the same reason. Again, we can construct a harmonic Green's kernel on ${\cal{K}}_k^-$. To do so we modify the kernel formula for the $k$-cylinders, in the way that we now insert a minus sign in the $k$-th coordinate of the expressions appearing in the sum. More precisely, 
\begin{theorem} Let $k < n-2$. Then the harmonic Green's kernel of harmonic sections on ${\cal{K}}_k^-$ with values in the trivial pin bundle can be expressed by
$$
H^{*}_k(x',y')= p_*\Big(\sum\limits_{m_1\underline{\omega}_1+...+m_k\underline{\omega}_k \in \Omega_k}  \Big|(\underline{x}-\underline{y}+\sum_{i=1}^{k-1} m_i\underline{\omega}_i + (-1)^{m_k} m_k\underline{\omega}_k)$$
$$+(x_{k+1}-y_{k+1}) e_{k+1} + ... + (x_{n-1}-y_{n-1}) e_{n-1} + (x_n-y_n) e_n  \Big|^{2-n}\Bigg).
$$
where $p_*$ again stands for the projection that maps a point $x \in \mathbb{R}^n$  to $x'=x \mod \sim^*$ which belongs to the manifold ${\cal{K}}_k^-$.  
\end{theorem}
The remaining cases $k\ge n-2$ can again be treated by re-grouping the terms of the series in a way coupling certain groups of periods together, which provides further degrees of convergence. This can be done in the same way as we did for the $k$-cylinders and the $n$-torus in \cite{KraRyan3}.


\begin{thebibliography}{999}
\bibitem{a} L. V. Ahlfors, {\it{M\"{o}bius transformations in $\mathbb{R}^{n}$ expressed through $2\times 2$ matrices of Clifford numbers}}, Complex Variables, 5, 1986, 215--224.
\bibitem{AH} P. Alexandroff and H. Hopf, {\it Topologie I}. Chelsea, Bronx, New York, 1935.
\bibitem{BBS} K. Becker, M. Becker, J. Schwarz, {\it String theory and $m$-theory}, Cambridge University Press, New York, 2007.
\bibitem{bds} F. Brackx, R. Delanghe and F. Sommen, {\it{Clifford Analysis}}, Pitman 76, London, 1982.
\bibitem{BCKR} E. Bulla, D. Constales, R.S. Krau{\ss}har and J. Ryan, {\it Dirac Type Operators for Arithmetic Subgroups of
Generalized Modular Groups}, J. reine angew. Math. 643, 2010, 1--19.
\bibitem{ca} D. Calderbank, {\it{Dirac operators and Clifford analysis on manifolds with boundary}}, Max Planck Institute for Mathematics, Bonn, preprint number 96-131, 1996.
\bibitem{cn} J. Cnops, {\it{An Introduction to Dirac Operators on Manifolds}}, Progress in Mathematical Physics, Birkh\"{a}user, Boston, 2002.
\bibitem{ER} M. Eastwood and J. Ryan, {\it Aspects of Dirac operators in analysis}, Milan J. Math. 75, 2007, 91--116.
\bibitem{f} T. Friedrich, {\it{Zur Abhangigheit des Dirac-operators von der Spin-Struktur}}, Colloq. Math., 48, 1984, 57--62.
\bibitem{GHS} K. G\"urlebeck, K. Habetha and W. Spr\"o{\ss}ig, {\it Holomorphic functions in the plane and $n$-dimensional space},
Birkh\"auser Verlag, Basel, 2008.
\bibitem{GS2} K.~G\"urlebeck and W.~Spr\"o{\ss}ig. {\it Quaternionic and Clifford Calculus for
Physicists and Engineers}, John Wiley \& Sons, Chichester-New York, 1997.
\bibitem{KraHabil} R. S. Krau{\ss}har, {\it Generalized analytic automorphic forms in hypercomplex spaces}, Frontiers in  Mathematics, Birkh\"auser, Basel, 2004. 
\bibitem{KraRyan1} R. S. Krau{\ss}har and John Ryan, {\it Clifford and harmonic analysis on cylinders and tori}, Rev. Matematica Iberoamericana, 21, 2005, 87--110.
\bibitem{KraRyan3} R. S.~Krau{\ss}har, Yuying Qiao and J. Ryan, {\it Harmonic, monogenic and hypermonogenic functions
on some conformally flat manifolds in $\mathbb{R}^n$ arising from special arithmetic groups of the Vahlen group}, Contemporary Mathematics, 370, 2005, 159--173.
\bibitem{KraRyan2} R.S. Krau{\ss}har and J. Ryan, {\it Some conformally flat spin manifolds, Dirac operators and Automorphic forms}, J. Math. Anal. Appl. 325 (1), 2007, 359--376. 
 \bibitem{lm} H. B. Lawson, Jr and M.-L. Michelson, {\it{Spin Geometry}}, Princeton University Press, Princeton, NJ, 1989.
\bibitem{Kuiper} N.H. Kuiper, {\it{On conformally flat spaces in the large}}, Ann. Math., (2) 50, 1949, 916--924.
\bibitem{lmq} C. Li, A. McIntosh and T. Qian, {\it{Clifford algebras, Fourier transforms, and singular convolution operators on Lipschitz surfaces}}, Revista Mathematica Iberoamericana, 10, 1994, 665--721.
 \bibitem{lr} H. Liu and J. Ryan, {\it{Clifford analysis techniques for spherical pde}} Journal of Fourier Analysis and its Applications, 8 No.~6,  2002, 535--563.
\bibitem{ma} M. Markel, {\it{Regular functions over conformal quaternionic manifolds}}, Commentationes Mathematicae Universitatis Carolinae, 22, 1981, 579--583.
 \bibitem{m} A. McIntosh, {\it{Clifford algebras, Fourier theory, singular integrals, and harmonic functions on Lipschitz domains}}, Clifford Algebras in Analysis and Related Topics, edited by J. Ryan, CRC Press, Boca Raton, 1996, 33--87.
 \bibitem{mp} R. Miatello and R. Podesta, {\it{Spin structures and spectra of $\mathbb{Z}_{2}$ manifolds}}, Math. Z., 247, 2004, 319--335.
 \bibitem{pp} P. Petersen, {\it{Riemannian Geometry}}, Graduate Texts in Mathematics, 171, Springer Verlag, New York, 1997.
\bibitem{pf} F. Pf\"{a}ffle, {\it{The Dirac spectrum of Bieberbach manifolds}}, J. Geom. Phys., 35, 2000, 367--385.
\bibitem{p} I. Porteous, {\it{Clifford Algebras and Classical Groups}}, Cambridge University Press, Cambridge, 1995.
\bibitem{r85} J. Ryan, {\it{Conformal Clifford manifolds arising in Clifford analysis}}, Proc. R. Ir. Acad., Sect. A 85,  1985, 1--23.
\bibitem{r2} J. Ryan, {\it{Conformally covariant operators in Clifford analysis}}, Zeitschrift f\"{u}r Analysis und ihre Anwendungen, 14, 1995, 677--704. 
\bibitem{Ry2003} J. Ryan, \textit{Cauchy kernels for some conformally flat manifolds}, Advances in analysis and geometry, Trends in math. Birkh\"auser, Basel, 2004, 149--160.
\bibitem{v} P. Van Lancker, {\it Clifford analysis on the sphere}, Clifford Algebras and their Applications in Mathematical Physics, edited by V. Dietrich et al, Kluwer, Dordrecht, 1998, 201--215.
\end{thebibliography}
\end{document}